\documentclass[graybox]{svmult}

\usepackage{helvet}         
\usepackage{courier}        
\usepackage{type1cm}        
\usepackage{makeidx}         
\usepackage{graphicx}        
\usepackage{multicol}        
\usepackage[bottom]{footmisc}
\usepackage{epstopdf}

\graphicspath{ 
{images/}
}
\usepackage{mathtools}

\makeindex             

\usepackage{graphics}

\usepackage{caption}
\usepackage{subcaption}
\usepackage{amsmath}
\usepackage{amssymb}
\usepackage{amsfonts}		
\usepackage{algorithm}
\usepackage[noend]{algpseudocode}

\newcommand{\cB}{\ensuremath{\mathcal{B}}}

\newcommand{\cD}{\ensuremath{\mathcal{D}}}

\newcommand{\cJ}{\ensuremath{\mathcal{J}}}
\newcommand{\cK}{\ensuremath{\mathcal{K}}}

\newcommand{\cM}{\ensuremath{\mathcal{M}}}
\newcommand{\cN}{\ensuremath{\mathcal{N}}}

\newcommand{\cS}{\ensuremath{\mathcal{S}}}

\newcommand{\cV}{\ensuremath{\mathcal{V}}}

\newcommand{\cX}{\ensuremath{\mathcal{X}}}

\newcommand{\bE}{\ensuremath{\mathbb{E}}}

\newcommand{\bR}{\ensuremath{\mathbb{R}}}

\DeclareMathOperator*{\argmax}{arg\,max}
\DeclareMathOperator*{\argmin}{arg\,min}

\DeclareMathOperator*{\vspan}{span}

\newcommand{\vphi}{\ensuremath{\phi}}

\newcommand{\test}{\ensuremath{\text{test}}}
\newcommand{\refl}{\ensuremath{\text{refl}}}
\newcommand{\core}{\ensuremath{\text{core}}}
\newcommand{\keff}{\ensuremath{k_{\text{eff}}}}

\newcommand{\eps}{\ensuremath{\varepsilon}}
\newcommand{\measure}{\sigma}
\newcommand{\solve}{\ensuremath{\textsc{solve}}}
\newcommand{\dist}{\operatorname{dist}}

\begin{document}

\title*{Stabilization of (G)EIM in presence of measurement noise: application to nuclear reactor physics}
\titlerunning{Stabilization of (G)EIM in presence of measurement noise}
\author{J.P.~Argaud, B.~Bouriquet, H.~Gong, Y.~Maday and O.~Mula}
\institute{J.P.~Argaud, B.~Bouriquet and H.~Gong \at \'{E}lectricit\'{e} de France, R\&D, 7 boulevard Gaspard Monge, 91120 Palaiseau, France \email{jean-philippe.argaud@edf.fr; bertrand.bouriquet@edf.fr; helin.gong@edf.fr}
\and Y.~Maday \at  Sorbonne Universit\'{e}, UPMC Univ Paris 06, UMR 7598, Labo. J.-L. Lions, 75005, Paris, France; Institut Universitaire de France; Div. of Applied Maths, Brown Univ., Providence R.I., USA \email{maday@ann.jussieu.fr}
\and O.~Mula \at Universit\'e Paris-Dauphine, PSL Research University, CNRS, UMR 7534, CEREMADE, 75016 Paris, France \email{mula@ceremade.dauphine.fr}
}


\maketitle

\abstract{The  Empirical Interpolation Method (EIM) and its generalized version (GEIM) can be used to approximate a physical system by combining data measured from the system itself and a reduced model representing the underlying physics. In presence of noise, the good properties of the approach are blurred in the sense that the approximation error no longer converges but even diverges. We propose to address this issue by a least-squares projection with constrains involving a some a priori knowledge of the geometry of the manifold formed by all the possible physical states of the system. The efficiency of the approach, which we will call Constrained Stabilized GEIM (CS-GEIM), is illustrated by numerical experiments dealing with the reconstruction of the neutron flux in nuclear reactors. A theoretical justification of the procedure will be presented in future works.}

\section{General overview and contribution of the paper}
For the sake of clarity, we shall start by formulating the goal of the paper in general terms containing  statements that will be clarified in the forthcoming sections.

Let $\cX$ be a Banach space over a domain $\Omega\subset \bR^d\ (d\geq 1)$ being equipped with the norm $\|.\|_\cX$. Our goal is to approximate functions $f$ from a given compact set $\cS\subset \cX$ which represent the states of a physical system taking place in $\Omega$. For this, any $f\in\cS$ will be approximated by combining two ingredients. The first is the use of a certain amount $m$ of measurements of $f$ collected directly from the system itself. We represent them as linear functionals of $\cX'$ (the dual of $\cX$) evaluated on $f$. The second ingredient is the use of a (family of) subspace(s) $V_n$ of finite dimension $n$ which is assumed to approximate well the set $\cS$. To limit the complexity in the approximation and also economize in the amount of sensors to place in the system, a desirable feature is to find appropriate sensors and the appropriate spaces $V_n$ for which $m$ and $n$ are moderate. A necessary hypothesis to allow this is to assume some properties on the geometry of $\cS$ expressed in terms of  a rapid decay  of the Kolmogorov $n$-width of $\cS$ in $\cX$,
\begin{equation*}
d_n(\cS,\cX) \coloneqq \inf_{ \substack{X\subset \cX \\ \dim(X)\leq n} } \max_{u\in\cS} \min_{v\in X} \Vert u - v \Vert_\cX .
\end{equation*}
Under this hypothesis on the decay of $d_n(\cS,\cX) $, one can in principle build a sequence $\{X_n\}_n$  s.t.
$\dist(\cS,X_n) \coloneqq  \max_{u\in\cS} \min_{v\in X} \Vert u - v \Vert_\cX \leq \varepsilon,
$
where $\dim(X_n) = n \equiv n(\varepsilon)$ is moderate.

Algorithms to build $\{X_n\}_n$ (or at least the first spaces in the sequence allowing to approximate beyond a given accuracy) and find appropriate linear functionals have been proposed in the community of reduced modeling (see \cite{MM2013,MPPY2014,MMT2016}). Note that, even if this is not required in the previous statements, the construction of the spaces $X_n$ is then recursive, i.e.  $X_{n-1}Ñ\subset X_n$. There, the approximation of $f\in\cS$ is done by interpolation or related approximations. The methods (in practice mainly based on a greedy procedure) are however not robust with respect to noise in the measurements and this paper introduces a constrained least squares approximation for which numerical experiments indicate its potential to address this obstruction.

\section{Mathematical setting}
Let us assume that $\cM = \overline{\hbox{span} (\cS)}$ (where the $\overline \cB $   denotes the closure in $\cX$ of the set $\cB$) admits a Schauder basis $\{q_i\}_i$, i.e., for every $f\in\cX$ there exists a unique sequence $\{c_i(f)\}$ of scalars such that
$\lim_{n\to\infty} \Vert f-\sum_{i=1}^n c_i(f) q_i \Vert_\cX = 0.
$
 For every $n\geq 1$, we define the $n$-dimensional subspace
$
X_n \coloneqq \vspan\{q_1,\dots,q_n\}.
$
Let us formulate in a different manner the hypothesis made involving the Kolmogorov $n$-width of $\cS$ in $\cX$: let us assume that the error in approximating the functions of $\cS$ in $X_n$ is
$
\max_{f\in\cS} \dist(f,X_n)\leq \eps_n,
$
where the sequence $(\eps_n)_n$ decays at a nice rate with $n$.

Let now $\{\lambda_i\}$ be the set of linear functionals of $\cX'$ (of unity norm in $\cX'$) such that for every $n\geq 1$, $\{\lambda_1,\dots,\lambda_n\}$ and $\{q_1,\dots,q_n\}$ are such that, for every $n\geq 1$ and every $1\leq j \leq n$,
$
\forall i, \  1\leq j \leq i \quad \lambda_j (q_i) = \delta_{i,j}.
$
For any $n\geq 1$, we can now define a (generalized) interpolation operator $\cJ_n:\cX\to X_n$ such that for all $f\in\cX$
$
\lambda_i(f) = \lambda_i\left( \cJ_n(f) \right),\quad i\in\{1,\dots,n\}.
$
By construction, for any $n\geq1$ and any $f\in\cX$,
$
\cJ_n(f)=\cJ_{n-1}(f)+c_n(f) q_n.
$
where
$
c_n(f) = \lambda_n\left( f-\cJ_{n-1}(f) \right)
$
and, for notational coherence, we set $\cJ_0=0$.

Using $\cJ_n$ to approximate the functions of $\cS$ yields the error bound
\begin{equation}
\max_{f\in\cS} \Vert f - \cJ_n(f) \Vert_\cX \leq \left( 1+ \Lambda_n \right) \eps_n,
\end{equation}
where
$
\Lambda_n \coloneqq \sup_{f\in\cX} \Vert \cJ_n(f) \Vert_\cX / \Vert f \Vert_\cX
$
is the Lebesgue constant. The value of $\Lambda_n$ diverges at a certain rate so the behavior of $\max_{f\in\cS} \Vert f - \cJ_n[f] \Vert$ with the dimension is dictated by the trade-off between the rate of divergence of $(\Lambda_n)$ (that is generally slow) and the convergence of $(\eps_n)$ (that is generally very fast). 
Also, for any $f \in \cS$
\begin{equation}\label{eq:1}
| c_n(f) | \leq (1+\Lambda_{n-1}) \eps_{n-1},\quad n\geq 1
\end{equation}
where $\Lambda_{0}=0$ and $\eps_0 =  \max_{f\in \mathcal S} \Vert f \Vert$.

For any $\alpha >0$, let us define the cone
$
\cK_n(\alpha) \coloneqq \{ v\in V_n\ :\ v=\sum_{i=1}^n c_i q_i\,\ |c_i|\leq \alpha (1+\Lambda_{i-1})\eps_{i-1} \}.
$
We have for any $n\geq 1$ and any  $f\in \cS$
$
\cJ_n(f) \in \cK_n(1)
$.
In presence of noise in the measurements, we assume that we receive values $\eta_1(f),\dots,\eta_n(f)$ such that $\eta_i(f) \sim \lambda_i(f)+\cN(0,\sigma^2)$ for $i\in\{1,\dots,n\}$. Interpolating from these values yields an element in $X_n$ denoted as $\cJ_n(f; \cN)$ that satisfies blurred error bound with respect to (\ref{eq:1}) that, depending on the precise definition of $\sigma$ 
and the norm of $\cX$ can be
\begin{equation}
\bE\left( \max_{f\in\cS} \Vert f - \cJ_n(f; \cN) \Vert \right) \leq \left( 1+ \Lambda_n \right) \eps_n + \left( 1+ \Lambda_n \right) \sqrt{n} 
\sigma.
\end{equation}

The second term of the bound diverges as $n$ increases and shows that the method is not asymptotically robust in presence of noise. An illustration of this can be found in the numerical results below. This motivates the search for other methods which would ideally yield a bound of the form $\left( 1+ \Lambda_n \right) \eps_n + \sigma$ and for which the error is asymptotically at the level of the noise  $\sigma$.

We propose to correct the interpolation operator by using more the structure of the manifold $\cS$ that, at the discrete level, is expressed in the fact that the approximation should belong to $\cK_n$. Indeed, the belonging of $\cJ_n(f; \cN)$ to $\cK_n$ is not satisfied any more except if there exists $\tilde f$ in $\cS$ such that $  \lambda_i(\tilde f) = \eta_i(f) $ for any $i, 1\le i \le n$ (which is rarely the case). In addition, in order to minimize the effect of the noise, we can increase the number of measurements and use $m$ larger
 than $n$ linear functional evaluations at a given dimension $n$. This leads to propose  a least-squares projection on $\cK_n$. For a given $n$, we now collect the values $\eta_{n,1}(f),\dots,\eta_{n,m(n)}(f)$ with $m(n)\geq n$ and such that
$
\eta_{n,i}(f) \sim \lambda_{i}(f)+\cN(0,\sigma^2),\quad 1\leq i\leq m(n).
$
Any $f\in\cS$ is now approximated by
\begin{equation}
A_n(f) = \argmin_{v\in\cK_n(\alpha)} \sum_{i=1}^{m(n)} \left( \lambda_{n,i}(v) - \eta_{n,i}(f) \right)^2
\label{eq:An}
\end{equation}
where $\alpha>1$ is suitabily chosen.

In this paper, we are running the greedy algorithm of the so-called Generalized Empirical Interpolation Method (GEIM, \cite{MM2013}), to generate a  basis $\{q_i\}_i$ and the linear functionals $\{\lambda_i\}_i$. This method is reported to have a nice behavior for the Lebesgue constant, at least in case it is trained on a set $\cS$ with small Kolmogorov dimension (see \cite{MMT2016}). This approach allows an empirical optimal selections of the positions of the sensors that provide (in case where no noise pollutes the measures)
a stable representation of the physical system. The precise algorithm is documented elsewhere (see  \cite{MM2013} and \cite{GABM2016}). Then, we approximate any $f\in\cS$ with the function $A_n(f)$ defined in \eqref{eq:An} with $\alpha = 2$. We call this scheme Constrained Stabilized GEIM (CS-GEIM). 

Note that the above approach could also be used with a POD approach to provide the imbedded spaces $\{X_n\}_n$ (that are more expensive to provide than the greedy GEIM approach but are more accurate) and well chosen linear functionals $\lambda_{n,i}$ the choice of which infer on the behavior of the Lebesgue constant $\Lambda_n$.

\section{Numerical results}
\label{sec:numerical}

\subsection{Modelling the physical problem}
\label{subsec::physical}

For the physical problem that we consider in this paper, the model is the two group neutron diffusion equation : the flux $\vphi$ has two energy groups $\vphi=(\vphi_1,\vphi_2)$. Index 1 denotes the high energy group and 2 the thermal energy one. These are modeled by the following parameter dependent PDE model :
\begin{equation}
\label{eq:diffusion}
\begin{cases}
&-\nabla\left(D_1\nabla\vphi_1\right)+(\Sigma_{a,1}+\Sigma_{s,1\to2})\vphi_1=\frac{1}{\keff} \left( \chi_1\nu\Sigma_{f,1}\vphi_1+ \chi_1\nu\Sigma_{f,2}\vphi_2 \right) \\
&-\nabla\left(D_2\nabla\vphi_2\right)+\Sigma_{a,2}\vphi_2-\Sigma_{s,1\to2}\vphi_1
=\frac{1}{\keff}\left( \chi_2\nu\Sigma_{f,1}\vphi_1+ \chi_2\nu\Sigma_{f,2}\vphi_2 \right),
\end{cases}
\end{equation}
here $\keff$ is the so-called multiplication factor and is not a data but an unknown of the problem \footnote{We omit here the technical details on the meaning of $\keff$ and refer to general references like \cite{Hebert2009}.}, 
and the given parameters are
\begin{itemize}
\item $D_i$ is the diffusion coefficient of group $i$ with $i\in\{1,2\}$.
\item $\Sigma_{a,i}$ is the macroscopic absorption cross section of group $i$.
\item $\Sigma_{s,1\to2}$ is the macroscopic scattering cross section from group 1 to 2.
\item $\Sigma_{f,i}$ is the macroscopic fission cross section of group $i$.
\item $\nu$ is the average number of neutrons emitted per fission.
\item $\chi_i$ is the fission spectrum of group $i$. 
\end{itemize}
they are condensed in
$
\mu
=
\{
D_1, D_2, \Sigma_{a,1}, \Sigma_{a,2}, \Sigma_{s,1\to2}, \nu\Sigma_{f,1}, \nu\Sigma_{f,2}, \chi_1, \chi_2
\}$.

 
 \noindent We assume that the parameters of our diffusion model range in, say, $$D_1 \in [D_{1,\min}, D_{1,\max}],\ D_2 \in [D_{2,\min}, D_{2,\max}],\dots,\chi_2\in[\chi_{2,\min}, \chi_{2,\max}],$$ then
$
\cD 
:=
[D_{1,\min}, D_{1,\max}]\times\dots\times[\chi_{2,\min}, \chi_{2,\max}]
$
is the set of all parameters and the set of all possible states of the flux is given by
\begin{equation}
\label{eq:manifold}
\cS
:=
\{
(\vphi_1,\vphi_2, P)(\mu)\ :\ \mu \in \cD
\},
\end{equation}
where the power $P(\mu)$ is defined from $(\vphi_1,\vphi_2)(\mu)$ as
$
P(\mu)(x) := \nu\Sigma_{f,1} \vphi_1(\mu)(x) + \nu\Sigma_{f,2} \vphi_2(\mu)(x),\quad \forall x\in\Omega
$
We assume (see \cite{CD2014} for elements sustaining this hypothesis) that the Kolmogorov-width decays rapidly, hence, it is possible to approximate all the  states of the flux (given by $\cS$) with an accuracy $\varepsilon$  in  well-chosen subspaces $X_n\subset \cX$ of relatively small dimension $n(\varepsilon)$.

To ensure enough stability in the reconstruction and minimize the approximation error, it is necessary to find the optimal placement of the sensors in the core. The selection is done with GEIM. If we denote
$
\measure (\vphi_i,x),\quad i\in \{1,2\},
$
the measurement of $\vphi_i$ at a position $x\in \Omega$ by a certain sensor, this measurement  can be modeled by a local average over $\vphi_i$ centered at $x\in \Omega$. Another possibility is to directly assume that the value $\vphi_i(x)$ at point $x$ is $\measure (\vphi_i,x)$\footnote{In the following part of this work, we directly assume that the value $\vphi_i(x)$ at point $x$ is $\measure (\vphi_i,x)$ as measurement.}. 
Note that, in principle, the measurement could depend on other parameters apart from the position. We could imagine for instance that we have sensors with different types of accuracy or different physical properties. This flexibility is not included in the current notation but   the reader will be able to extrapolate from the current explanations.

A specificity of the approach here is that $\cS$ is composed of vectorial functions $(\vphi_1,\vphi_2, P)(\mu)$. We deliberately choose to take measurements only on one of the components (say $\vphi_2$) and thus reconstruct the whole field $\vphi_1,\ \vphi_2$ and $P$ with the only knowledge of thermal flux measurements. 

Another specificity of our approach is on the spatial location of the measurements. We consider two cases:
\begin{itemize}
\item Case I: the sensors can be placed at any point in the domain of definition of $\vphi_2$.
\item Case II: the admissible sensor locations are restricted to be deployed in a restricted part of that domain;
\end{itemize}
We have already reported in \cite{ABGMM2016} that these two specificities are well supported by the (G)EIM approach, as long as the greedy method is taught to achieve the goal of reconstructing the whole field $(\vphi_1,\vphi_2, P)(\mu)$.

Our aim here is to show that the noise can be controlled through our Constrained Stabilized (G)EIM approach.

\subsection{Description of the PARCS 2D IAEA benchmark}
\label{subsec::benchmark}

\begin{figure}[!h]
\includegraphics[width=0.6\textwidth]{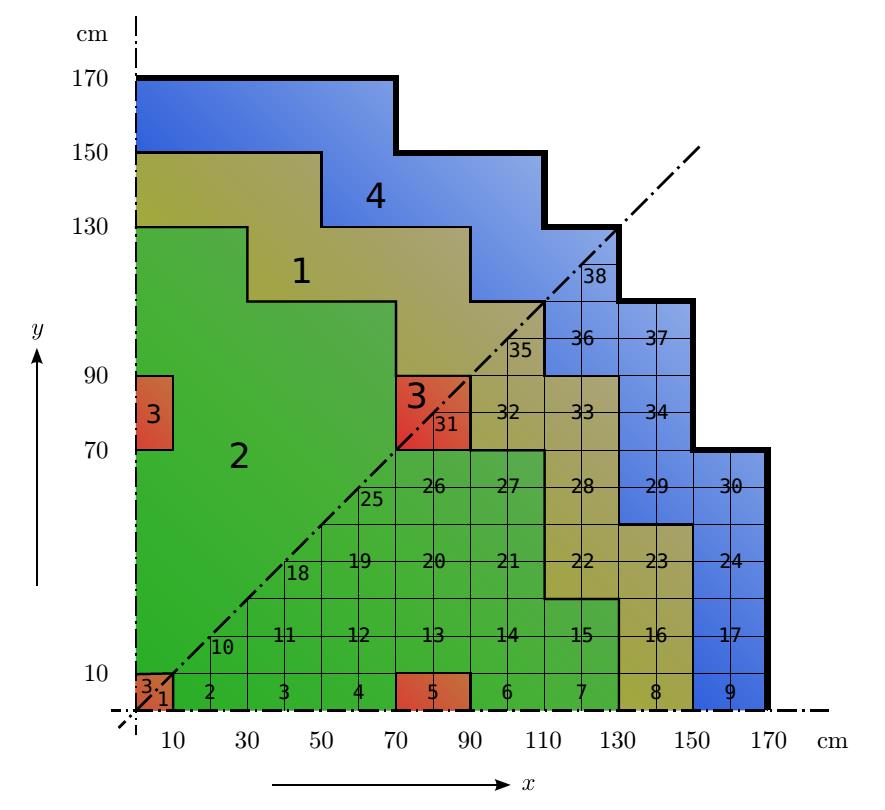}
\centering
\caption{Geometry of 2D IAEA benchmark, upper octant: region assignments, lower octant: fuel assembly identification (from reference \cite{2D-IAEA-Benchmark}). }
\label{fig::core configuration} 
\end{figure}

We consider the classical 2D IAEA Benchmark Problem \cite{Benchmark}, the core geometry which can be seen in figure \ref{fig::core configuration}. The problem conditions and the requested results are stated in page 437 of reference \cite{Benchmark}. It is identified with the code 11-A2, and its descriptive title is ¡°Two-dimensional LWR Problem¡±, also known as ¡°2D IAEA Benchmark Problem¡±. This problem represents the mid-plane $z=190~cm$ of the 3D IAEA Benchmark Problem, that is used by references \cite{PARCS} and show in application within \cite{2D-IAEA-Benchmark}.

The reactor domain is $\Omega=\text{region}(1,2,3,4)$. The core and the reflector are $\Omega_{\core}=\text{region}(1,2,3)$ and $\Omega_{\refl}=\text{region}(4)$ respectively. We consider only the value of $D_1|_{\Omega_{\refl}}$ in the reflector $\Omega_{\refl}$ as a parameter (so $p=1$ and $\mu=D_1|_{\Omega_{\refl}}$). We assume that $D_1|_{\Omega_{\refl}}\in[1.0, 3.0]$. The rest of the coefficients of the diffusion model \eqref{eq:diffusion} (including $D_1|_{\Omega_{\core}}$) are fixed to the values indicated on table \ref{tab:coefs2D}. In principle, one could also consider these coefficients as parameters but we have decided to focus only on $D_1|_{\Omega_{\refl}}$ because of its crucial role in the physical estate of the core: its variation can be understood as a change in the boundary conditions in $\Omega_{\core}$ which, up to a certain extent, allows to compensate the bias of the diffusion model with respect to reality. We shall report in a future paper more extended variations of the parameters.

\begin{table}
\caption{Coefficient values: macroscopic cross sections (units are not stated in the original reference, but they are assumed to be in $cm^{-1}$ or cm as appropriate.)}
\label{tab:coefs2D}       
%
%
\begin{tabular}{p{1cm}p{2.2cm}p{1.2cm}p{1.2cm}p{1.2cm}p{1.2cm}p{1.2cm}p{1.6cm}}
\hline\noalign{\smallskip}
Region & $D_1$ & $D_2$ & $\Sigma_{1\to 2}$ & $\Sigma_{a1}$ & $\Sigma_{a2}$ & $\nu \Sigma_{f2}$ & Material$^a$ \\
\noalign{\smallskip}\svhline\noalign{\smallskip}
1 & 1.5 & 0.4 & 0.02 & 0.01 & 0.080 & 0.135 & Fuel 1\\
2 & 1.5 & 0.4 & 0.02 & 0.01 & 0.085 & 0.135 & Fuel 2\\
3 & 1.5 & 0.4 & 0.02 & 0.01 & 0.130 & 0.135 & Fuel 2 + Rod\\
4 & [1.0, 3.0] or 2.0$^b$ & 0.3 & 0.04 & 0 & 0.010 & 0 & Reflector\\
\noalign{\smallskip}\hline\noalign{\smallskip}
\end{tabular}
$^a$ Axial buckling $B_{zg}^2=0.8\cdot10^{-4}$ for all regions and energy groups. \\
$^b$ here $2.0$ is the exact value from reference \cite{2D-IAEA-Benchmark}.
\end{table}

\subsection{Hypothesis of (CS-)GEIM application}
\label{subsec::application}

We propose to reconstruct $(\vphi_1,\vphi_2,P)$ as explained above (see \cite{ABGMM2016})
i.e., $\vphi_2$ will be approximated with its direct interpolant $\cJ_n[\vphi_2]$ while $\vphi_1$ and $P$ will be reconstructed from the measurements of $\vphi_2$, using the same coefficients in a coherent basis set. These are denoted as $\widetilde\cJ_n[\vphi_1]$ and $\widetilde\cJ_n[P]$.

%
%
%
%

\subsection{Numerical results}
\label{subsec::result}

\begin{figure}[H]
\centering
\begin{subfigure}[b]{5truecm}
\includegraphics[scale=0.40]{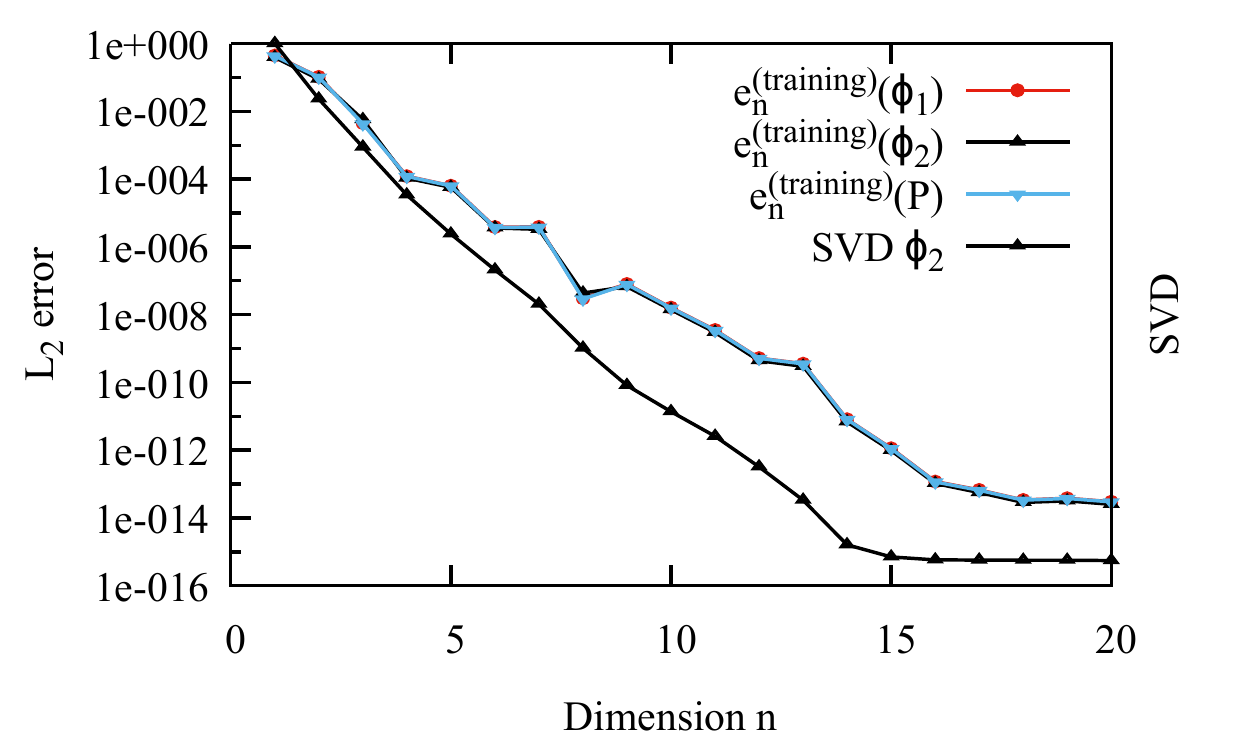}
\caption{Decay of SVD modes and of $e^{\text{(training)}}_n(\vphi_2)$.}
\label{fig:2d:svdVSgreedyL2}
\end{subfigure}%
~
\begin{subfigure}[b]{5truecm}
\includegraphics[scale=0.40]{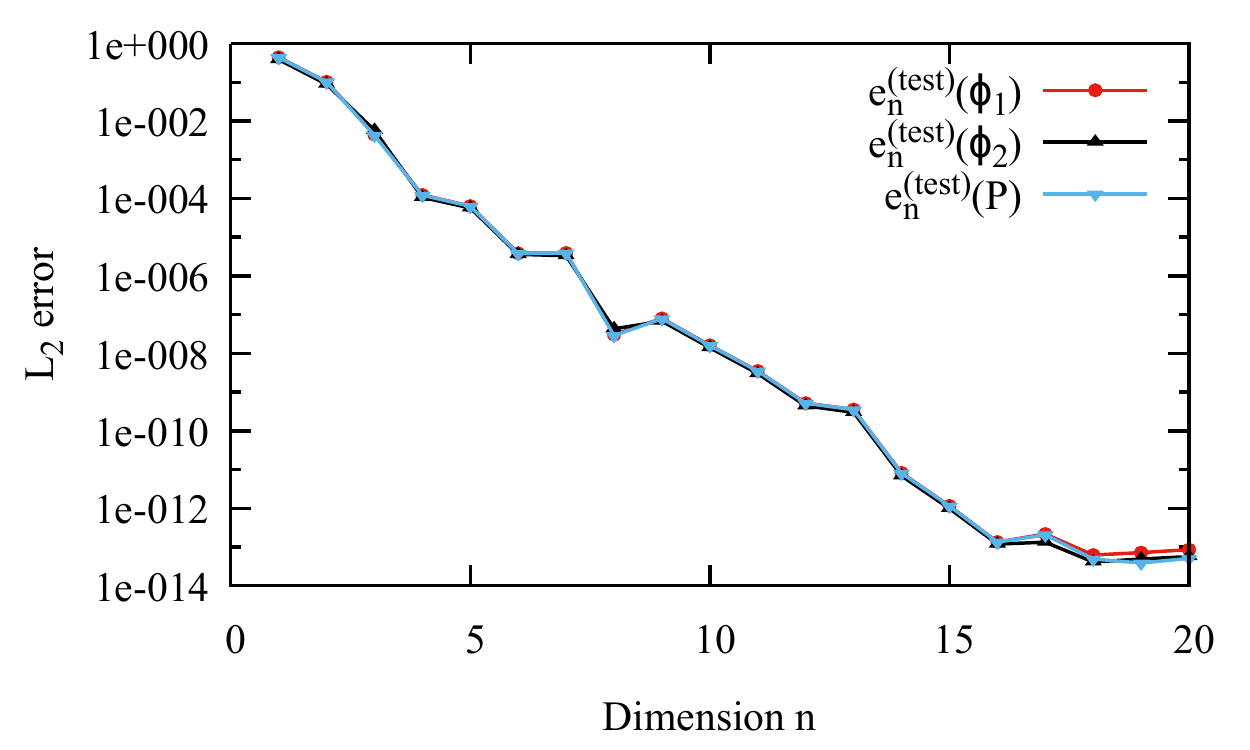}
\caption{Decay of $e^{\text{(test)}}_n(\vphi_1)$, $e^{\text{(test)}}_n(\vphi_2)$ and $e^{\text{(test)}}_n(P)$.}
\label{fig:2d:testSetL2}
\end{subfigure}
\caption{Case II, $L^2(\Omega)$ norm: Reconstruction of $(\vphi_1,\vphi_2,P)(\mu)$ with $\left(\widetilde \cJ_n[\vphi_1],\cJ_n[\vphi_2],\widetilde \cJ_n[P] \right)(\mu)$.}
\label{fig:2d:decayL2}
\end{figure}

Let us now turn to the analysis of the results. We study the performance of the reconstruction strategy by considering first of all the decay
of the errors
\begin{equation}
e^{\text{(training)}}_n(\vphi_2) 
:= \max_{\mu \in \cD^{\text{(training)}}} \Vert \vphi_2(\mu)-\cJ_n[\vphi_2](\mu) \Vert_{L^2(\Omega)}
\end{equation}
in the greedy algorithm. Since both case I and case II yield very similar results, we only present plots of case II for the sake of concision. In figure \ref{fig:2d:svdVSgreedyL2}, the decay is compared to an indicator of the optimal performance in $L^2(\Omega)$ which is obtained by a singular value decomposition of the snapshots $\vphi_2(\mu),\ \forall\mu\in\cD^{\text{(training)}}$. Note that $e^{\text{(training)}}_n(\vphi_2)$
decays at a similar rate as the SVD which suggests that GEIM behaves in a quasi-optimal way (see \cite{MMT2016}).
We now estimate the accuracy to reconstruct $(\vphi_1,\vphi_2,P)(D_1|_{\Omega_{\refl}})$ for any $D_1|_{\Omega_{\refl}}\in[0.5,2.0]$ which does not
necessary belong to the training set of snapshots. For this, we consider a test set of 300 parameters $\cD^{\text{(\test)}}$ different from $\cD^{\text{(training)}}$ and compute the errors
\begin{equation}
\label{eq:estimateError}
\begin{cases}
e^{\text{(test)}}_n(\vphi_1) 
&:= \max_{\mu \in \cD^{(\test)}} \Vert \vphi_1(\mu)-\widetilde \cJ_n[\vphi_1](\mu) \Vert_{L^2(\Omega)} \\
e^{\text{(test)}}_n(\vphi_2) 
&:= \max_{\mu \in \cD^{(\test)}} \Vert \vphi_1(\mu)-\cJ_n[\vphi_2](\mu) \Vert_{L^2(\Omega)} \\
e^{\text{(test)}}_n(P) 
&:= \max_{\mu \in \cD^{(\test)}} \Vert P(\mu)-\widetilde \cJ_n[P](\mu) \Vert_{L^2(\Omega)}
\end{cases}
\end{equation}
The decay of the errors \eqref{eq:estimateError} is plotted in figure \ref{fig:2d:testSetL2}. The fact that $e^{\text{(test)}}_n(\vphi_2)$ decays very similarly to $e^{\text{(training)}}_n(\vphi_2)$ confirms that the set of 300 training snapshots was representative enough of the whole manifold $\cS$. Also, the fast decay of $e^{\text{(test)}}_n(\vphi_1)$ and $e^{\text{(test)}}_n(P)$ shows that the use of the operator $\widetilde \cJ_n$ to approximate $\vphi_1$ and $P$ is accurate enough.

\begin{figure}[H]
\centering
\begin{subfigure}[b]{5truecm}
\includegraphics[scale=0.40]{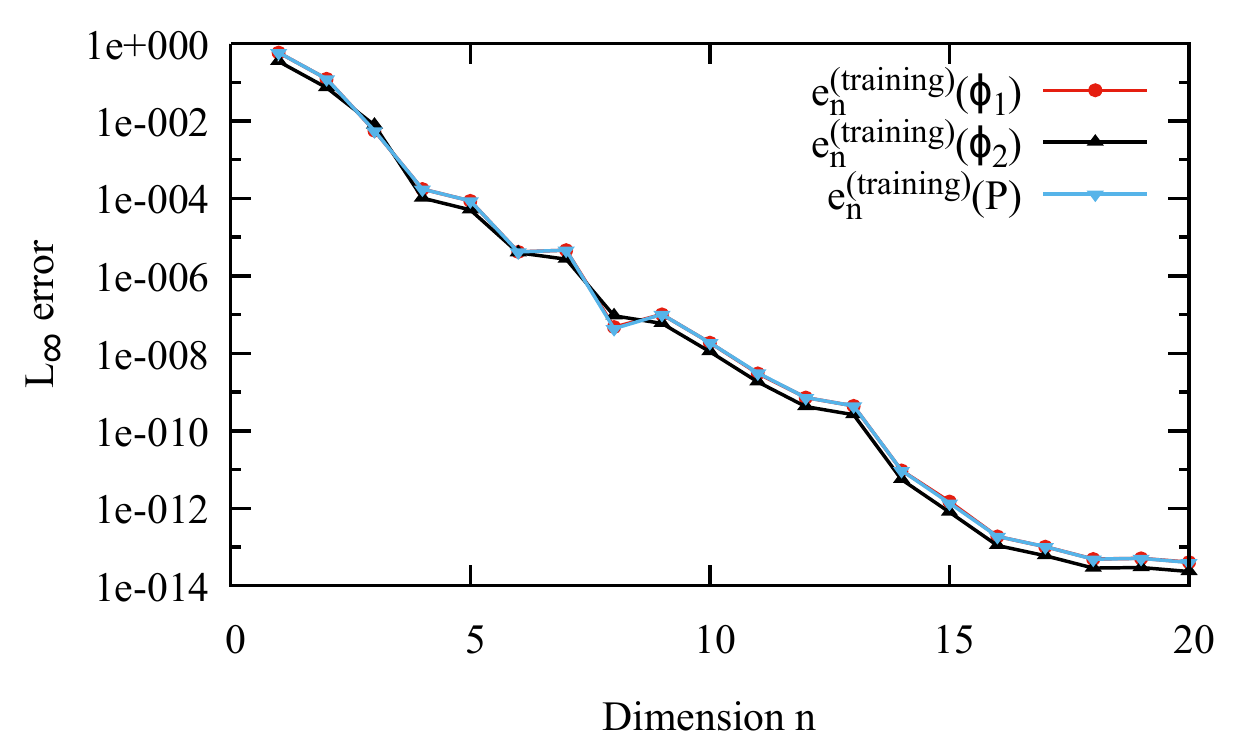}
\caption{Decay of the greedy errors $e^{\text{(training)}}_n(\vphi_2)$.}
\label{fig:2d:greedyLinf}
\end{subfigure}%
~
\begin{subfigure}[b]{5truecm}
\includegraphics[scale=0.40]{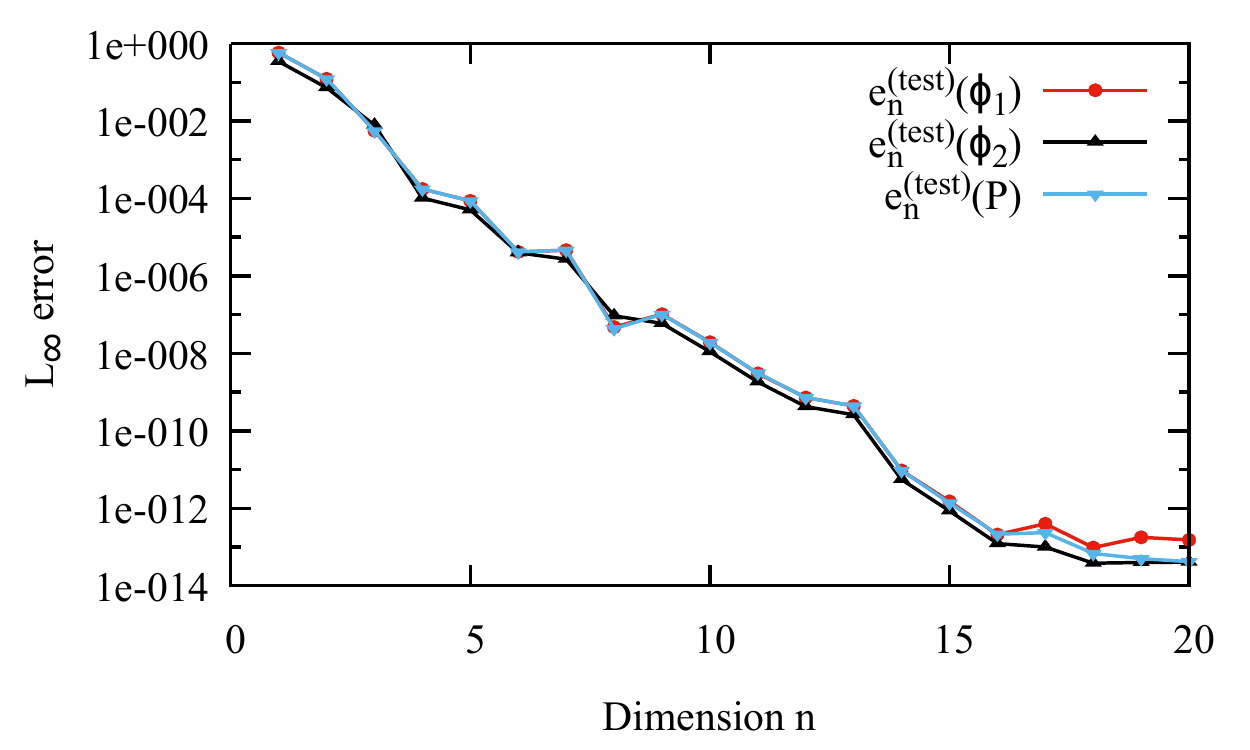}
\caption{Decay of $e^{\text{(test)}}_n(\vphi_1)$, $e^{\text{(test)}}_n(\vphi_2)$ and $e^{\text{(test)}}_n(P)$.}
\label{fig:2d:testSetLinf}
\end{subfigure}
\caption{Case II, $L^{\infty}(\Omega)$ norm: Reconstruction of $(\vphi_1,\vphi_2,P)(\mu)$ with $\left(\widetilde \cJ_n[\vphi_1],\cJ_n[\vphi_2],\widetilde \cJ_n[P] \right)(\mu)$.}
\label{fig:2d:decayLinf}
\end{figure}

\begin{figure}[H]
\centering
\begin{subfigure}[b]{5truecm}
\includegraphics[scale=0.40]{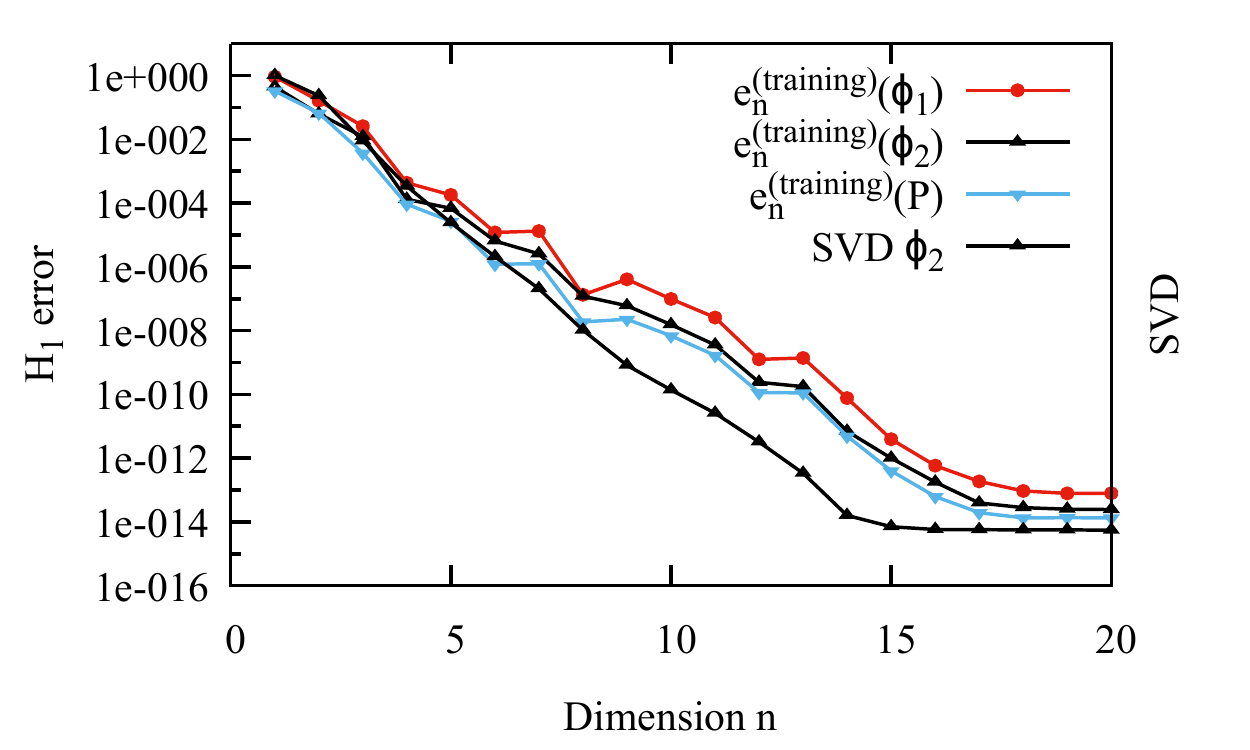}
\caption{Decay of SVD modes and of $e^{\text{(training)}}_n(\vphi_2)$.}
\label{fig:2d:svdVSgreedyH1}
\end{subfigure}%
~
\begin{subfigure}[b]{5truecm}
\includegraphics[scale=0.40]{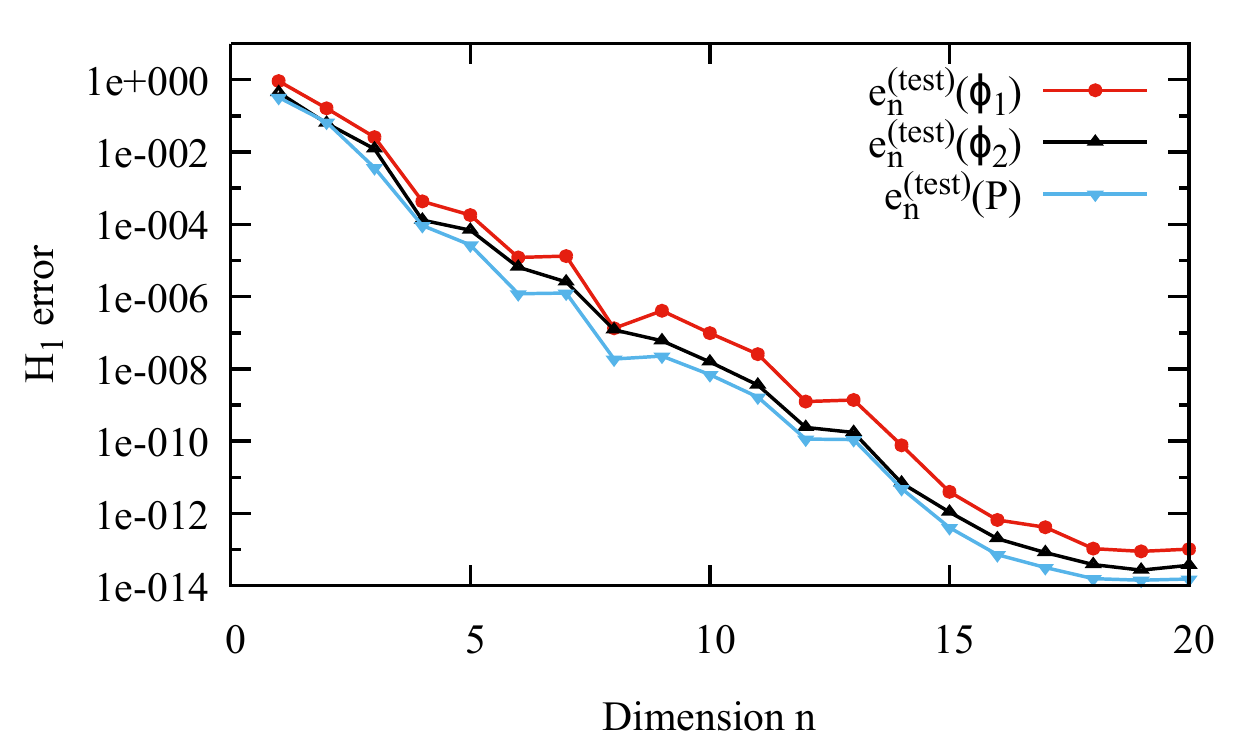}
\caption{Decay of $e^{\text{(test)}}_n(\vphi_1)$, $e^{\text{(test)}}_n(\vphi_2)$ and $e^{\text{(test)}}_n(P)$.}
\label{fig:2d:testSetH1}
\end{subfigure}
\caption{Case II, $H^1(\Omega)$ norm: Reconstruction of $(\vphi_1,\vphi_2,P)(\mu)$ with $\left(\widetilde \cJ_n[\vphi_1],\cJ_n[\vphi_2],\widetilde \cJ_n[P] \right)(\mu)$.}
\label{fig:2d:decayH1}
\end{figure}

Instead of working with $L^2(\Omega)$, it is also possible to work with other norms (provided some spacial regularity in the manifold). A particularly relevant case in neutronics is $L^\infty(\Omega)$. Figure \ref{fig:2d:decayLinf} shows the results of the reconstruction procedure when working in this norm and figure \ref{fig:2d:decayH1} shows the behavior when considering the  $H^1(\Omega)$ and working with its classical  semi-norm $|u|_{H^1(\Omega)}=\int_{\Omega} |\nabla u|^2$.



\begin{figure}[H]
\centering
\begin{subfigure}[b]{5truecm}
\includegraphics[scale=0.26]{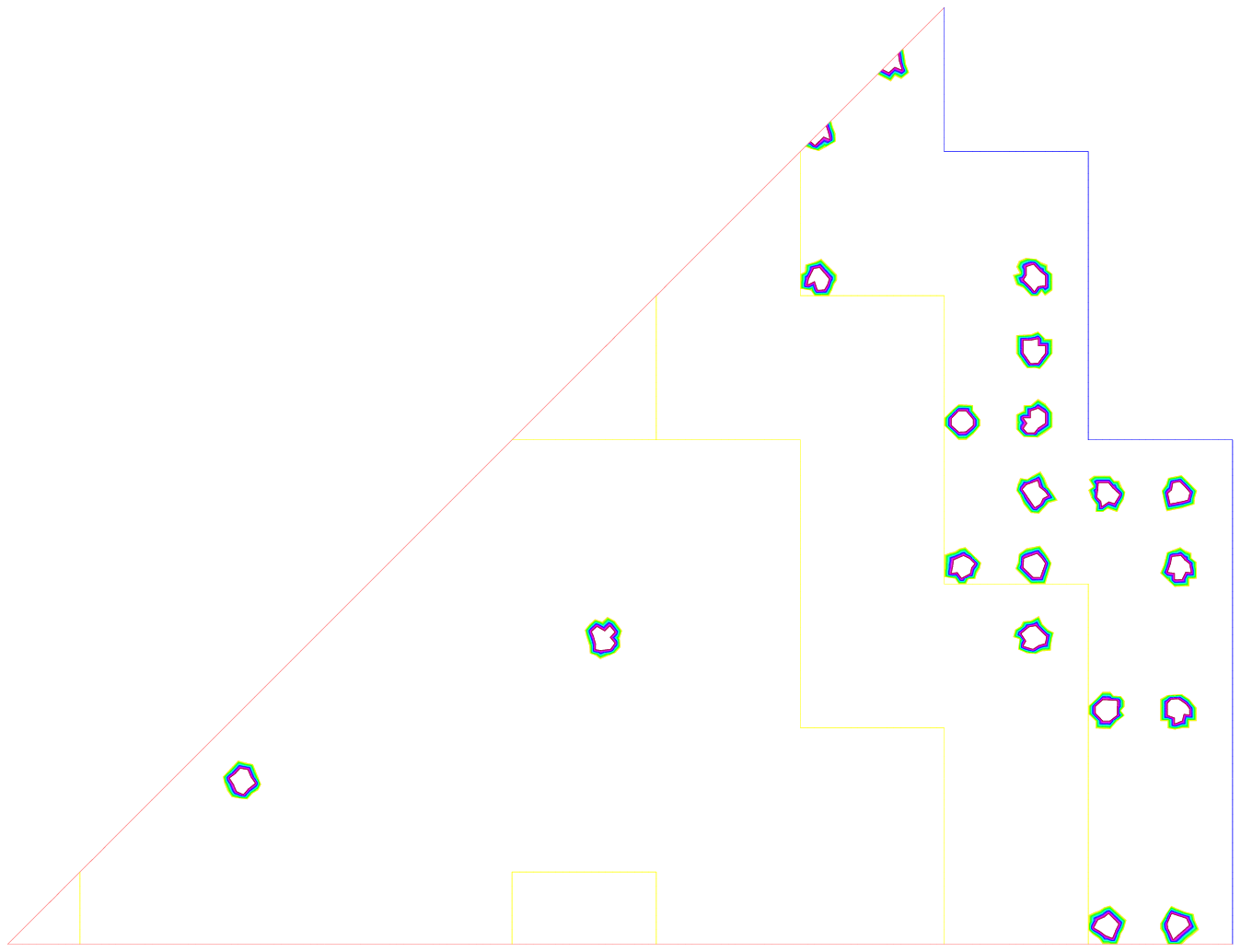}
\caption{Case I (selection in $\Omega$)}
\label{fig:pointsCaseI}
\end{subfigure}%
~
\begin{subfigure}[b]{5truecm}
\includegraphics[scale=0.26]{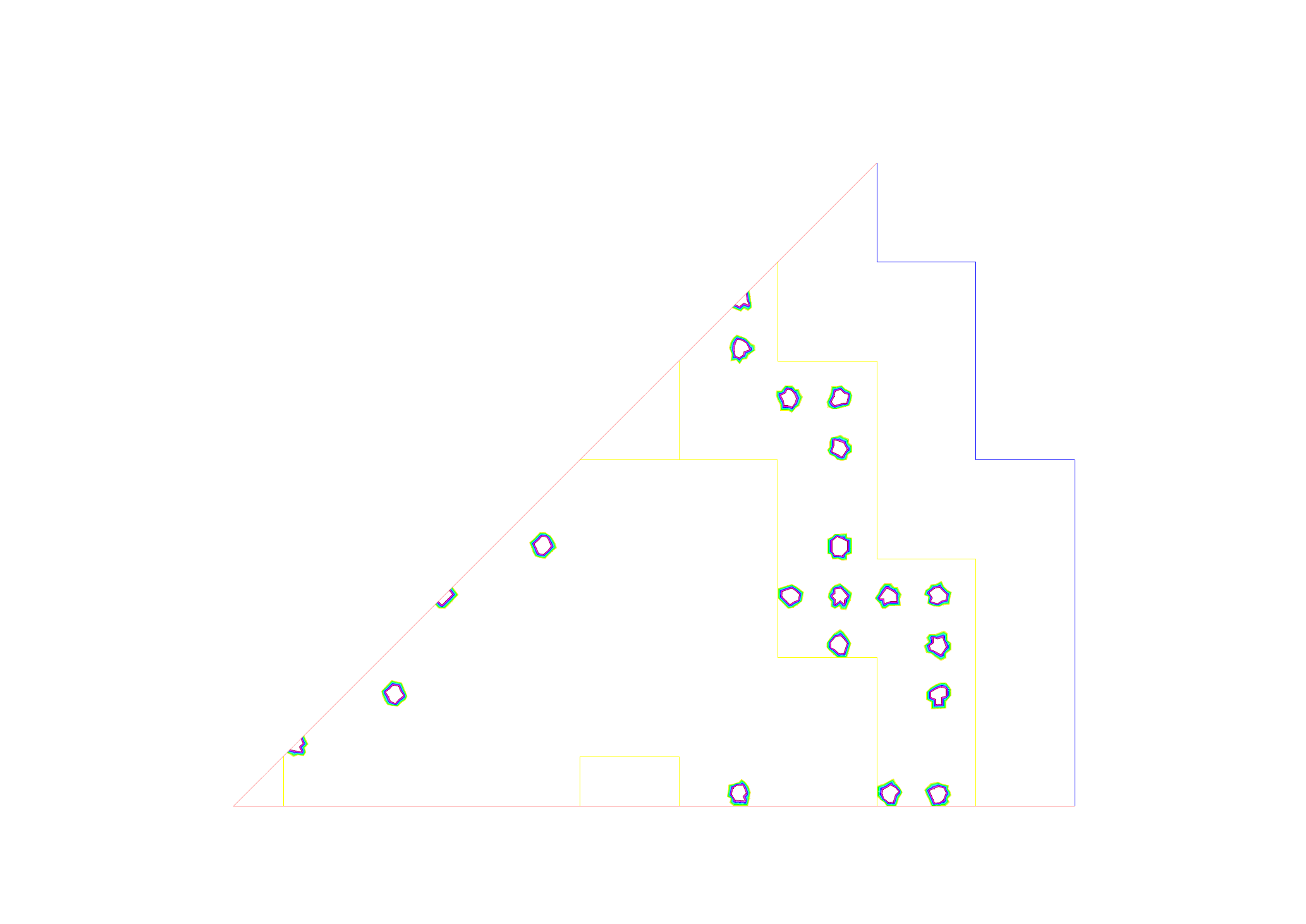}
\caption{Case II (selection in $\Omega_{\core}$)}
\label{fig:pointsCaseII}
\end{subfigure}
\caption{Locations of the sensors chosen by the greedy EIM algorithm.}
\label{fig::outputPointsGreedy}
\end{figure}


\begin{figure}[H]
\centering
\begin{subfigure}[b]{5truecm}
\includegraphics[scale=0.40]{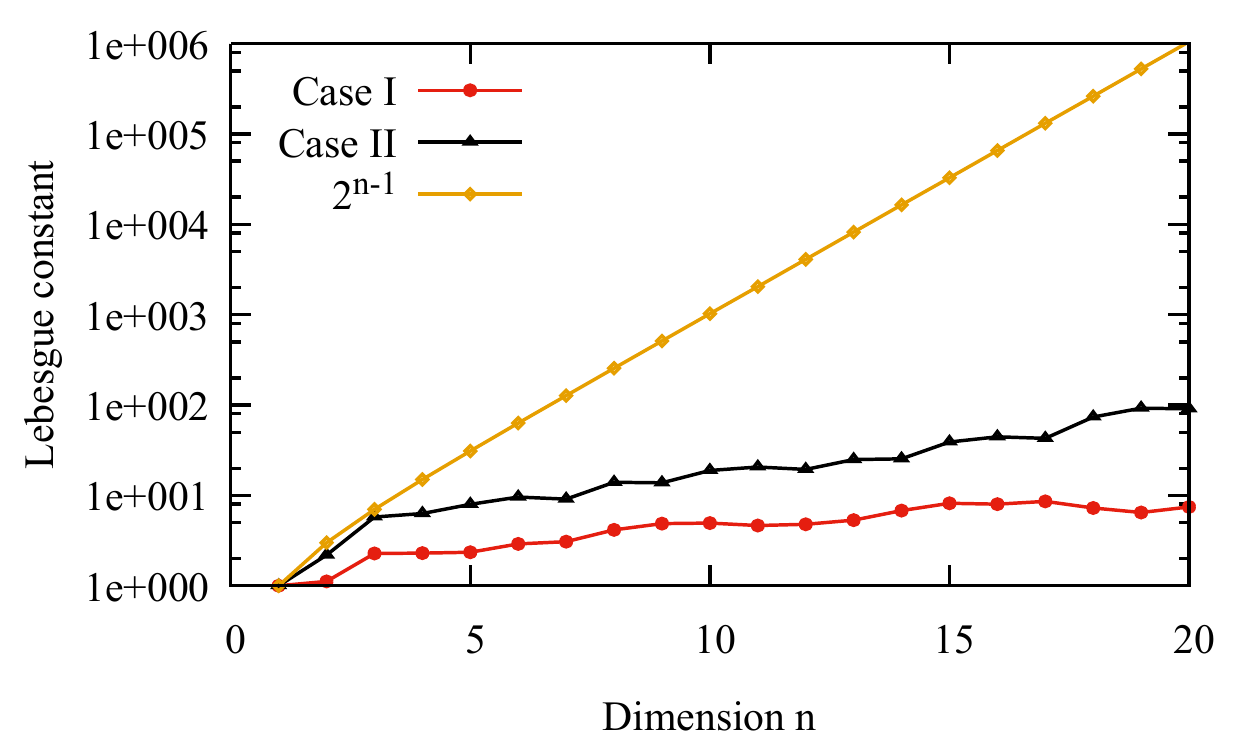}
\caption{The lebesgue constant for Case I and Case II}
\label{fig:2d:lebesgue}
\end{subfigure}%
~
\begin{subfigure}[b]{5truecm}
\includegraphics[scale=0.40]{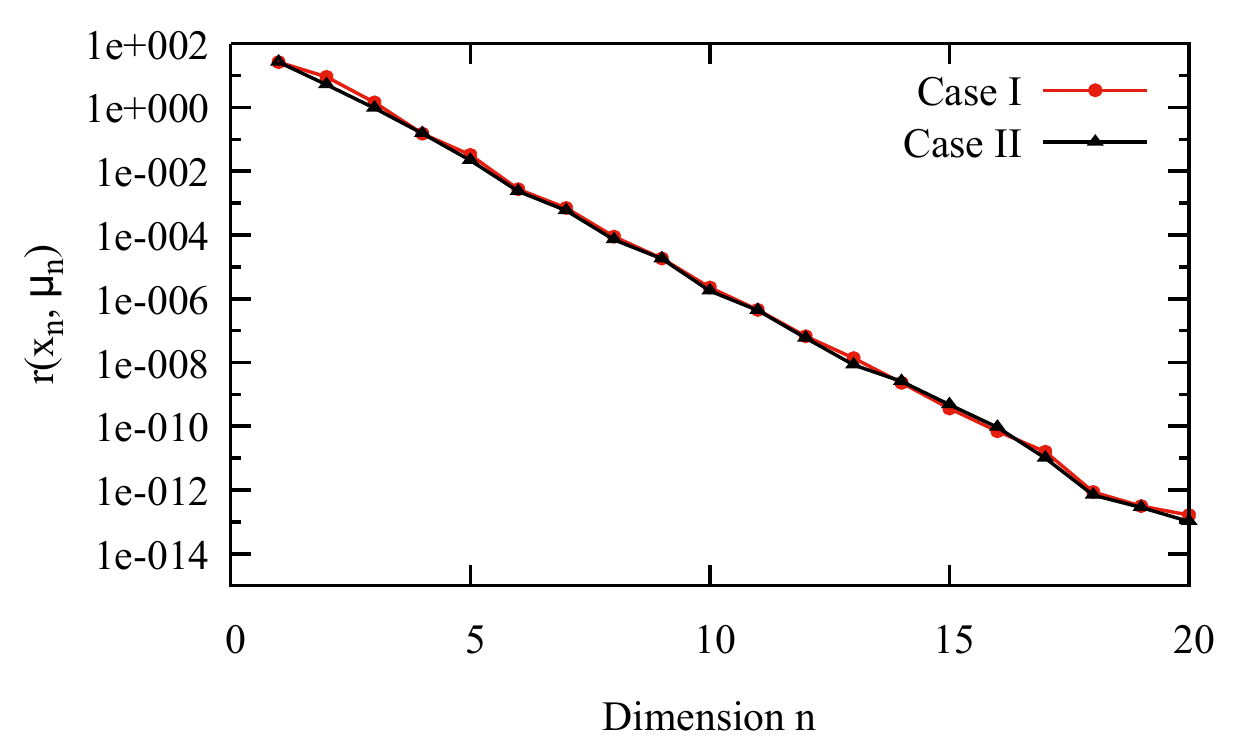}
\caption{The coefficients upper limits $r_n(x_n,\mu_n)$, for Case I and Case II.}
\label{fig:2d:coefficient}
\end{subfigure}
\caption{The lebesgue constant and $r_n(x_n,\mu_n)$ from GEIM}
\label{fig:2d:lebesgue-coefficient}
\end{figure}

Figure \ref{fig:2d:lebesgue} shows the behavior of the Lebesgue constants in both cases. We can find that 1) the Lebesgue constant increases with GEIM interpolation function dimension, 
2) if detectors are limited in a domain part (Case II), the Lebesgue constant gets worse, as an effect of the extrapolation that is required here, nevertheless the increase is still moderate. 


Figure \ref{fig:2d:coefficient} shows the coefficients upper limits described as $r_n(x_n,\mu_n)\equiv (1+\Lambda_{n-1}) \eps_{n-1}$ (see \eqref{eq:1}) (so  $|c_n| \le r_n(x_n,\mu_n) $) for Case I and Case II, which decreases quickly with $n$. 

\begin{figure}[H]
\centering
\begin{subfigure}[b]{5truecm}
\includegraphics[scale=0.40]{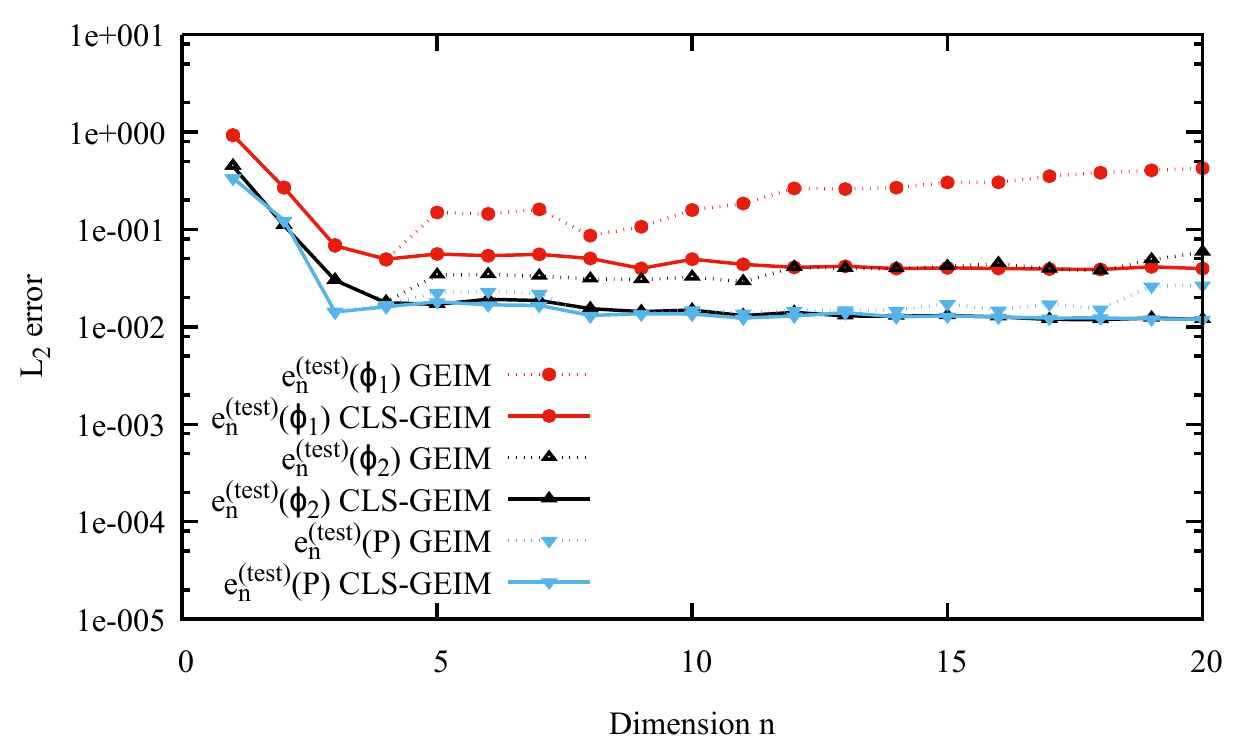}
\caption{Case I}
\label{fig:2d:1e-2casI}
\end{subfigure}%
~
\begin{subfigure}[b]{5truecm}
\includegraphics[scale=0.40]{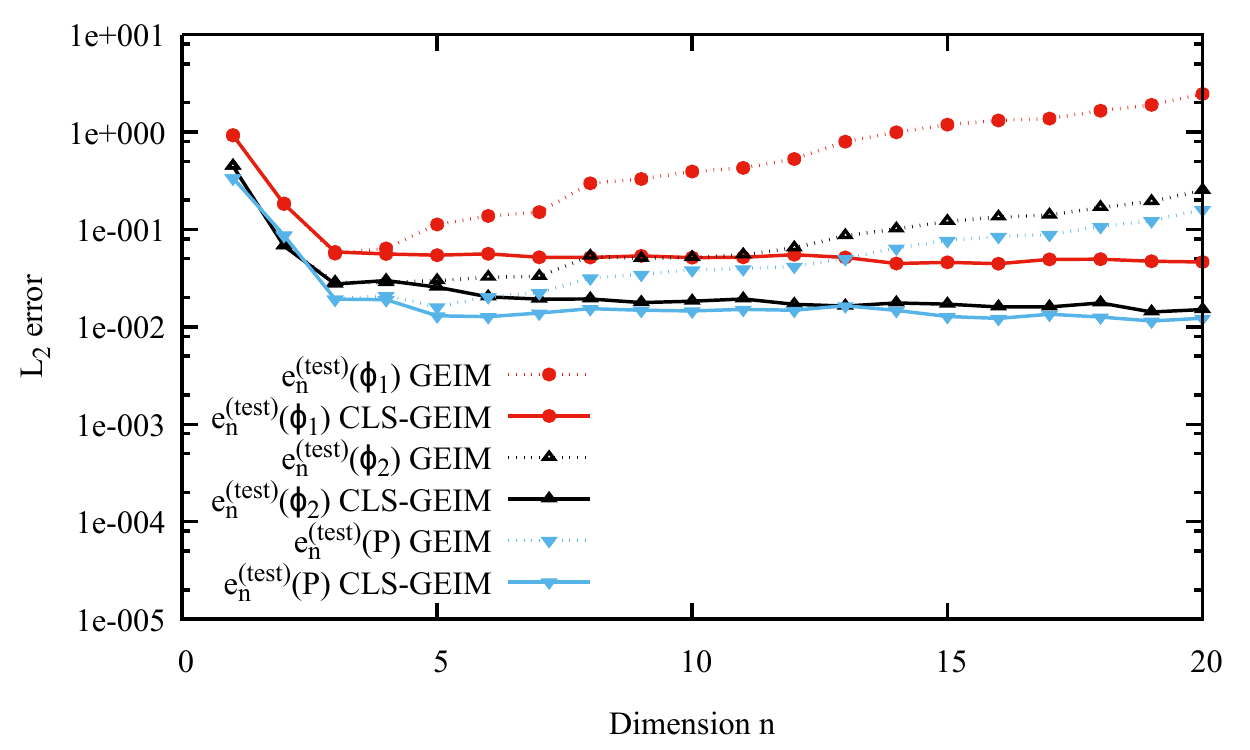}
\caption{Case II}
\label{fig:2d:1e-2casII}
\end{subfigure}
\caption{$L^2(\Omega)$ norm: decay of $e^{\text{(test)}}_n(\vphi_1)$, $e^{\text{(test)}}_n(\vphi_2)$ and $e^{\text{(test)}}_n(P)$. Reconstruction of $(\vphi_1,\vphi_2,P)(\mu)$ with $\left(\widetilde \cJ_n[\vphi_1],\cJ_n[\vphi_2],\widetilde \cJ_n[P] \right)(\mu)$, with noise amplitude $10^{-2}$.}
\label{fig:2d:1e-2}
\end{figure}
\begin{figure}[H]
\centering
\begin{subfigure}[b]{5truecm}
\includegraphics[scale=0.40]{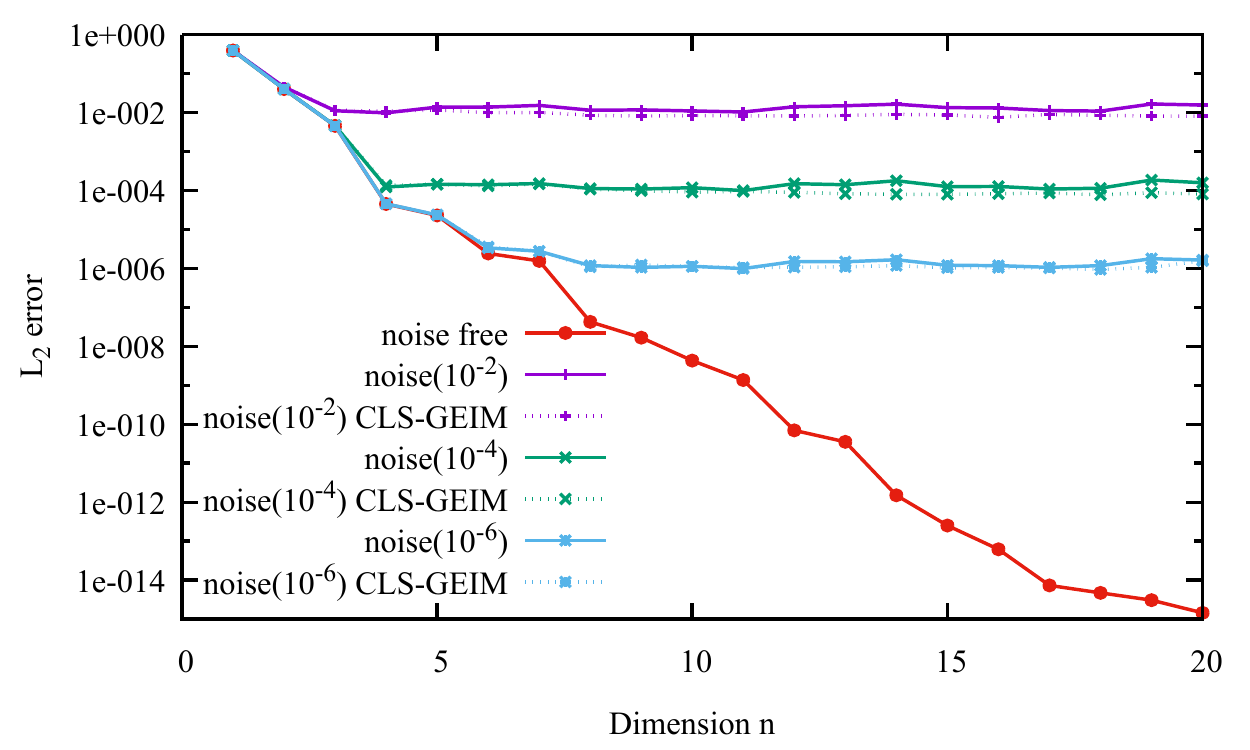}
\caption{Case I}
\label{fig:2d:CS-greedyL2-case-I}
\end{subfigure}%
~
\begin{subfigure}[b]{5truecm}
\includegraphics[scale=0.40]{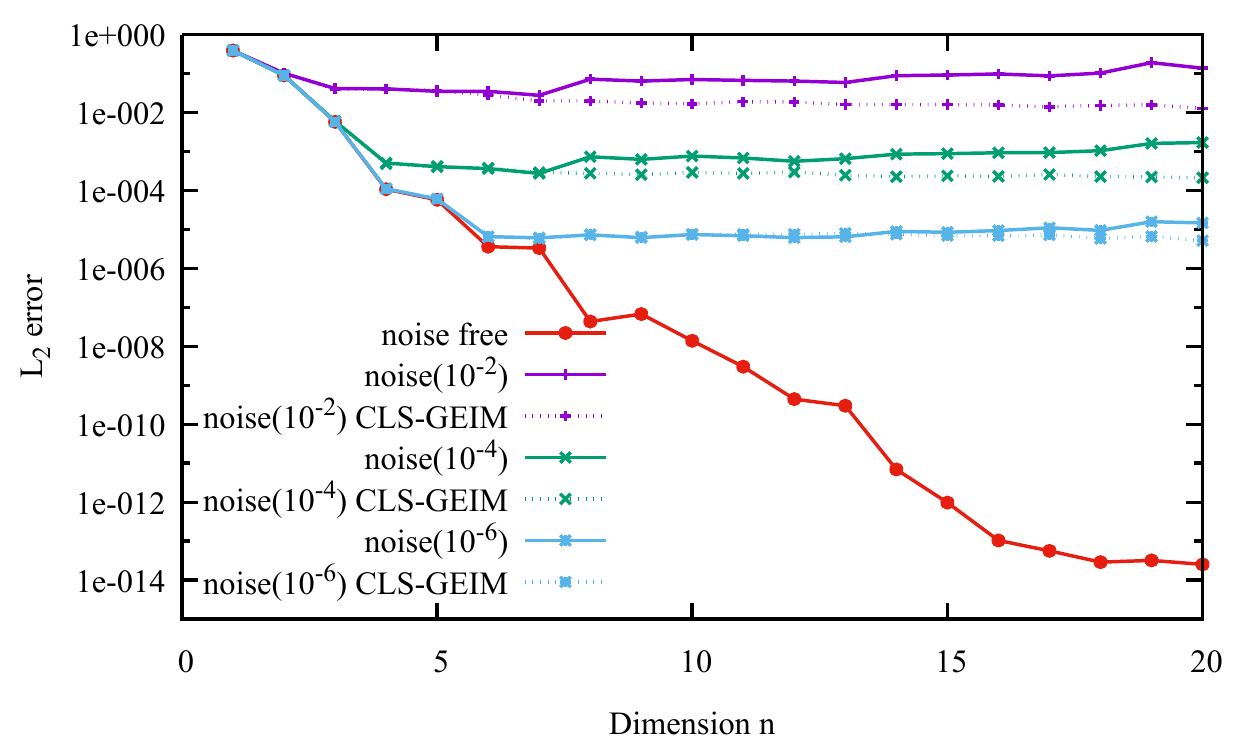}
\caption{Case II}
\label{fig:2d:CS-greedyL2-case-II}
\end{subfigure}
\caption{$L^2(\Omega)$ norm: decay of $e^{\text{(training)}}_n(\vphi_2)$, with noise amplitude $10^{-2}, 10^{-4}, 10^{-6}$.}
\label{fig:2d:CS-greedyL2}
\end{figure}
\begin{figure}[H]
\centering
\begin{subfigure}[b]{5truecm}
\includegraphics[scale=0.40]{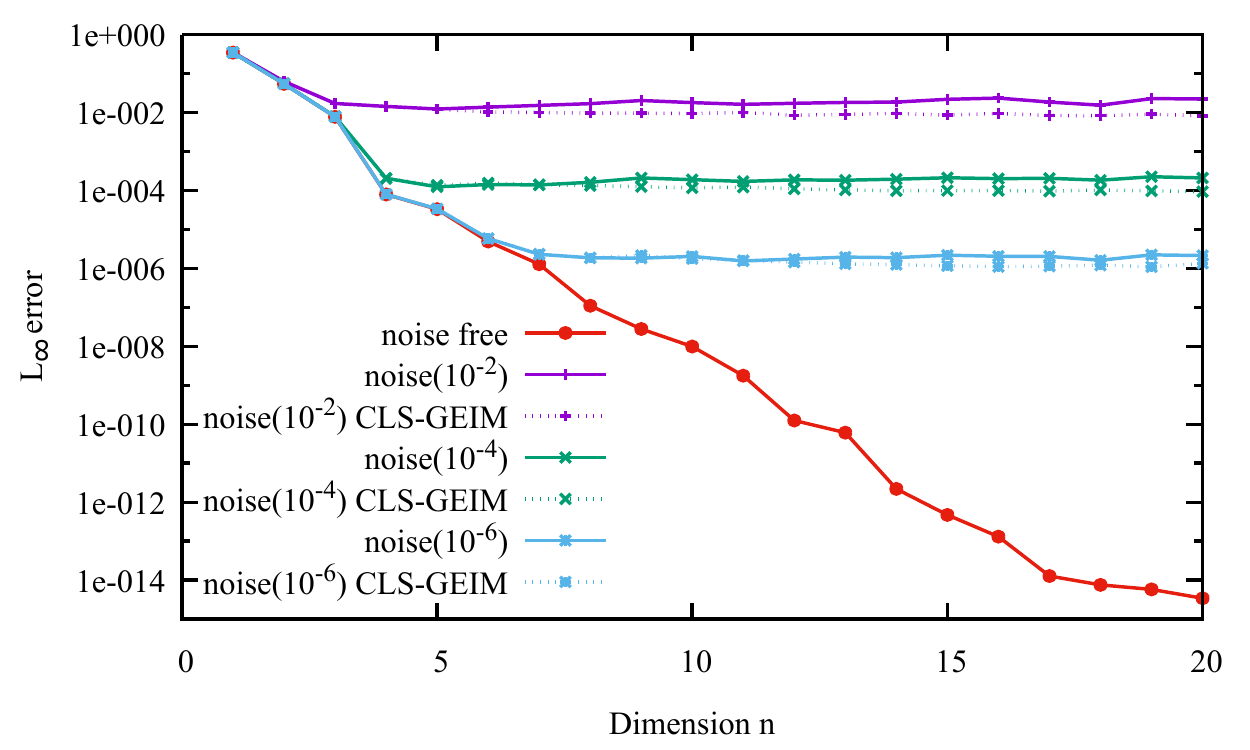}
\caption{Case I}
\label{fig:2d:CS-greedyLinf-case-I}
\end{subfigure}%
~
\begin{subfigure}[b]{5truecm}
\includegraphics[scale=0.40]{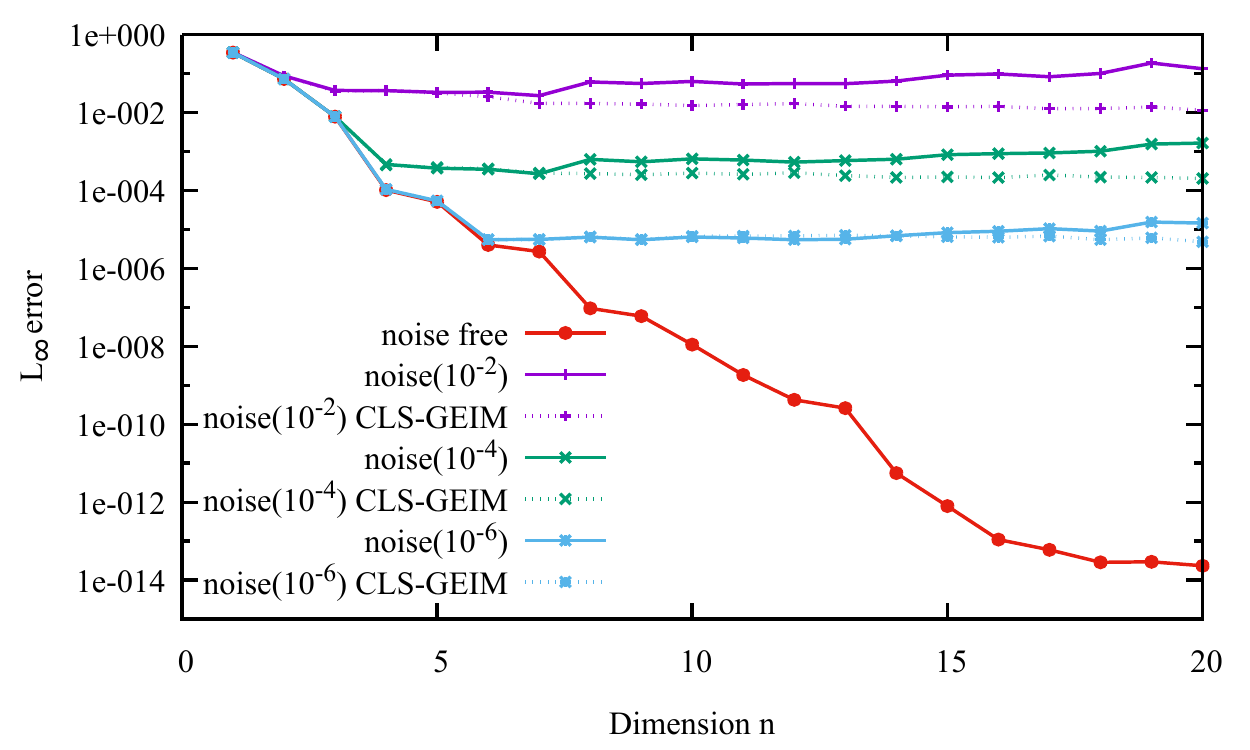}
\caption{Case II}
\label{fig:2d:CS-greedyLinf-case-II}
\end{subfigure}
\caption{$L^{\infty}(\Omega)$ norm: decay of $e^{\text{(training)}}_n(\vphi_2)$, with noise amplitude $10^{-2}, 10^{-4}, 10^{-6}$.}
\label{fig:2d:CS-greedyLinf}
\end{figure}
\begin{figure}[H]
\centering
\begin{subfigure}[b]{5truecm}
\includegraphics[scale=0.40]{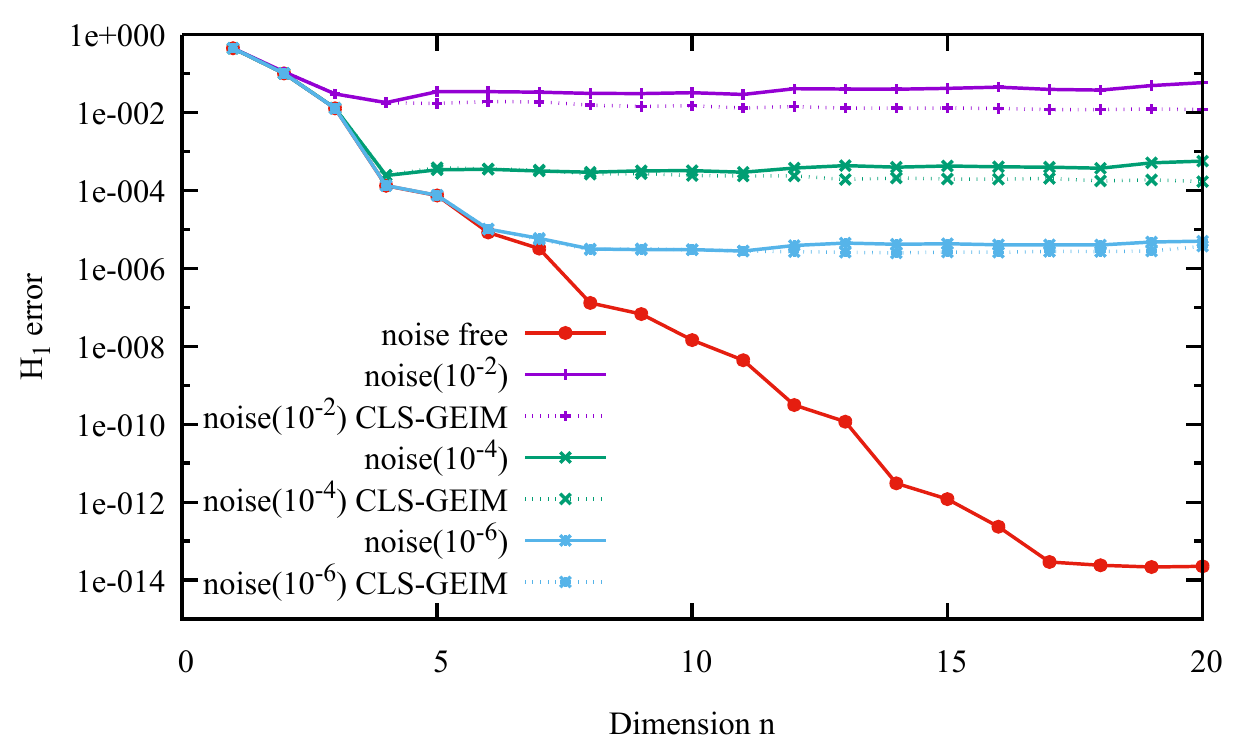}
\caption{Case I}
\label{fig:2d:CS-greedyH1-case-I}
\end{subfigure}%
~
\begin{subfigure}[b]{5truecm}
\includegraphics[scale=0.40]{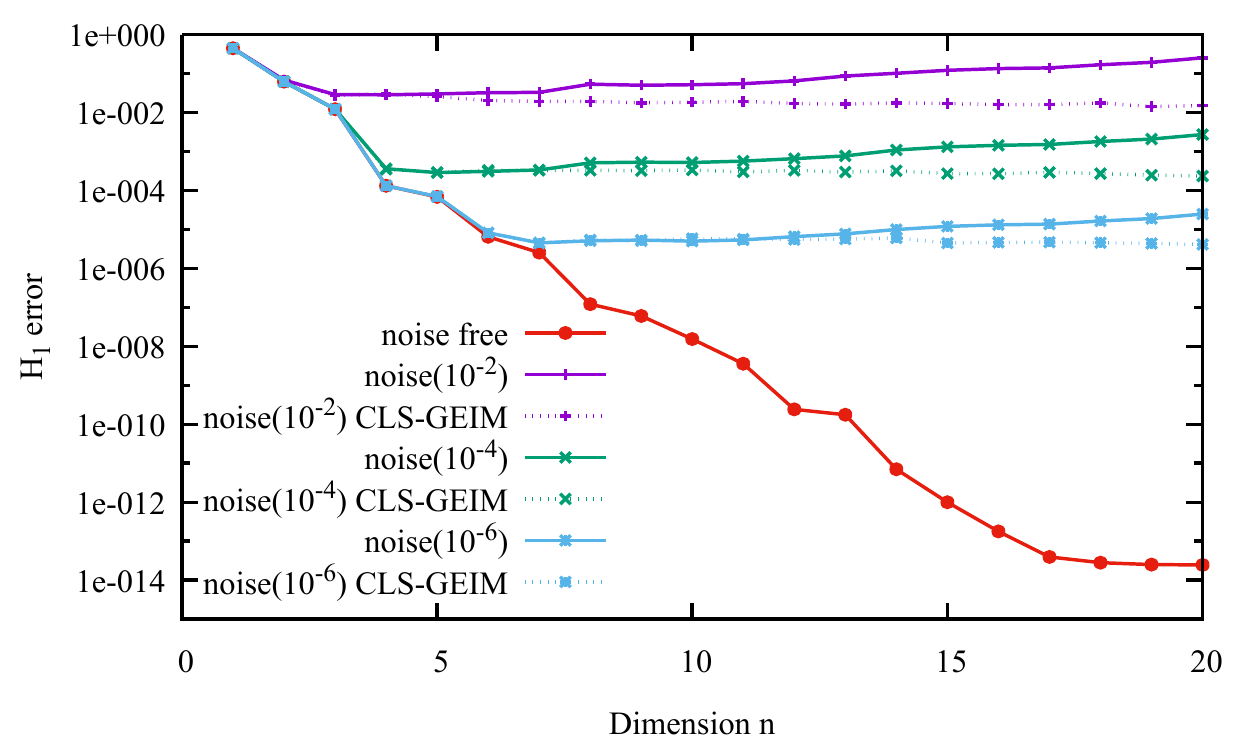}
\caption{Case II}
\label{fig:2d:CS-greedyH1-case-II}
\end{subfigure}
\caption{$H^1(\Omega)$ norm: decay of $e^{\text{(training)}}_n(\vphi_2)$, with noise amplitude $10^{-2}, 10^{-4}, 10^{-6}$.}
\label{fig:2d:CS-greedyH1}
\end{figure}

We still focus on the 300 parameters $\cD^{(test)}$, and compute the errors with equation (\ref{eq:estimateError}), for each test case, we perform the interpolation process with CS-GEIM a number of times.  
 Figure \ref{fig:2d:1e-2} shows the averaged $L^2(\Omega)$ norm for the decay of $e^{\text{(test)}}_n(\vphi_1)$, $e^{\text{(test)}}_n(\vphi_2)$ and $e^{\text{(test)}}_n(P)$,
 with noise amplitude $10^{-2}$ for Case I and Case II. For different measurement noise amplitude, the 
 averaged errors in $L^2(\Omega)$ norm, $L^{\infty}(\Omega)$ norm and $H^1(\Omega)$ norm are shown in  figure \ref{fig:2d:CS-greedyL2}, figure \ref{fig:2d:CS-greedyLinf} and figure \ref{fig:2d:CS-greedyH1}  respectively, for Case I and Case II. The main conclusions are: in the noisy case, i) CS-GEIM improves the interpolation, with the error comparable to the noise input level, ii) in extrapolation case, CS-GEIM reduces the interpolation error dramatically, which extends GEIM practical use.

\begin{figure}[H]
\centering
\includegraphics[scale=0.80]{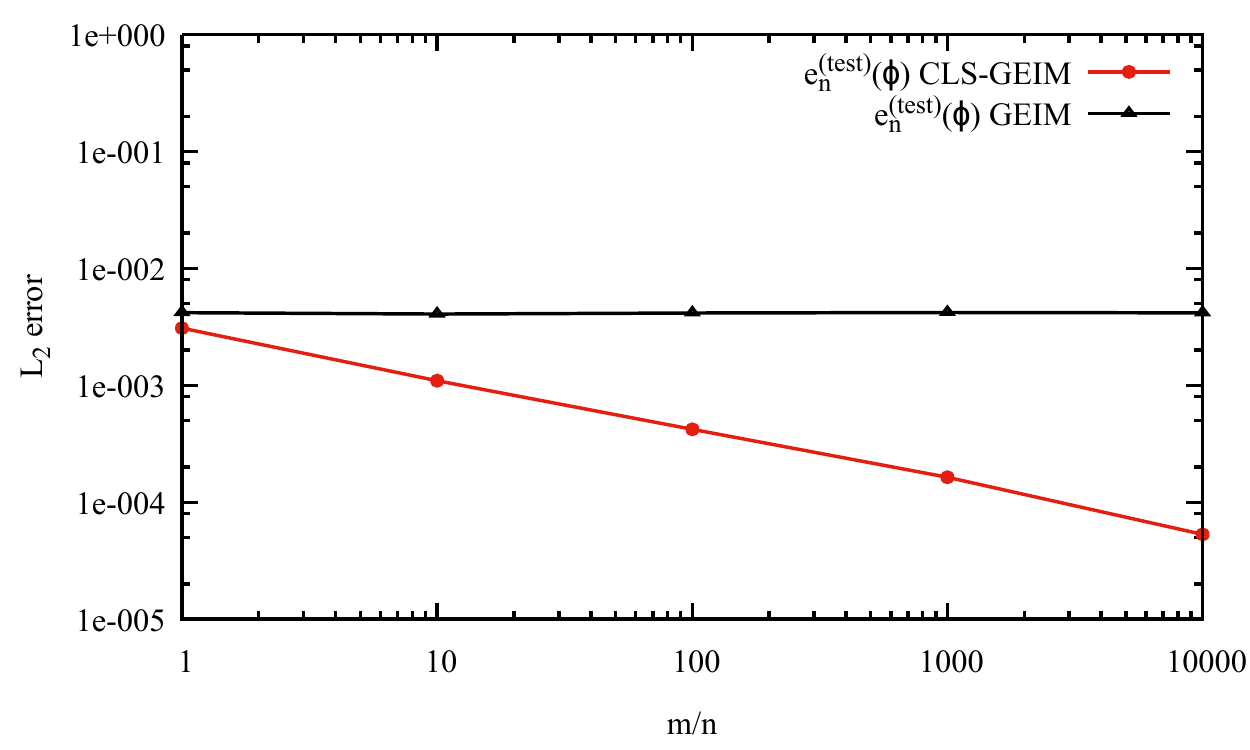}
\caption{CS-GEIM with different $n/m$ ratio, the input noise level is $10^{-2}$, with the function $g(x,\mu) \equiv  ((x_1-\mu_1)^2+(x_2-\mu_2)^2)^{-1/2}$, the error converges with $\sim n^{-\frac{1}{2}}$.}
\label{fig:2d:moremeasurement func2}
\end{figure}

If we take more measurements with fixed number of interpolation functions, the ratio $n/m$ of the number of measurements $n$ to the number of interpolation functions $m$ increases, so it is expected to have the same effect than to repeat independent measure at the same point in order to measure the evaluation of the measure. We consider the analytical function  $g(x,\mu) \equiv \cV((x_1,x_2);(\mu_1,\mu_2))\equiv ((x_1-\mu_1)^2+(x_2-\mu_2)^2)^{-1/2}$ for $x\in\Omega \equiv ]0,1[^2$ and $\mu\in\cD \equiv [-1,-0.01]^2$; we choose for $\cD^{(training)}$ a uniform discretization sample of 400 points\footnote{We replace the synthetic neutron problem here by the above analytical function so as to be able to have a more thorough and extensive numerical analysis}. Then we change the ratio $n/m$ of the number of measurements $n$ to the number of interpolation functions $m$ with CS-GEIM process, figure \ref{fig:2d:moremeasurement func2} also shows the error converges with $\sim n^{-\frac{1}{2}}$.

\begin{figure}[H]
\centering
\includegraphics[scale=0.80]{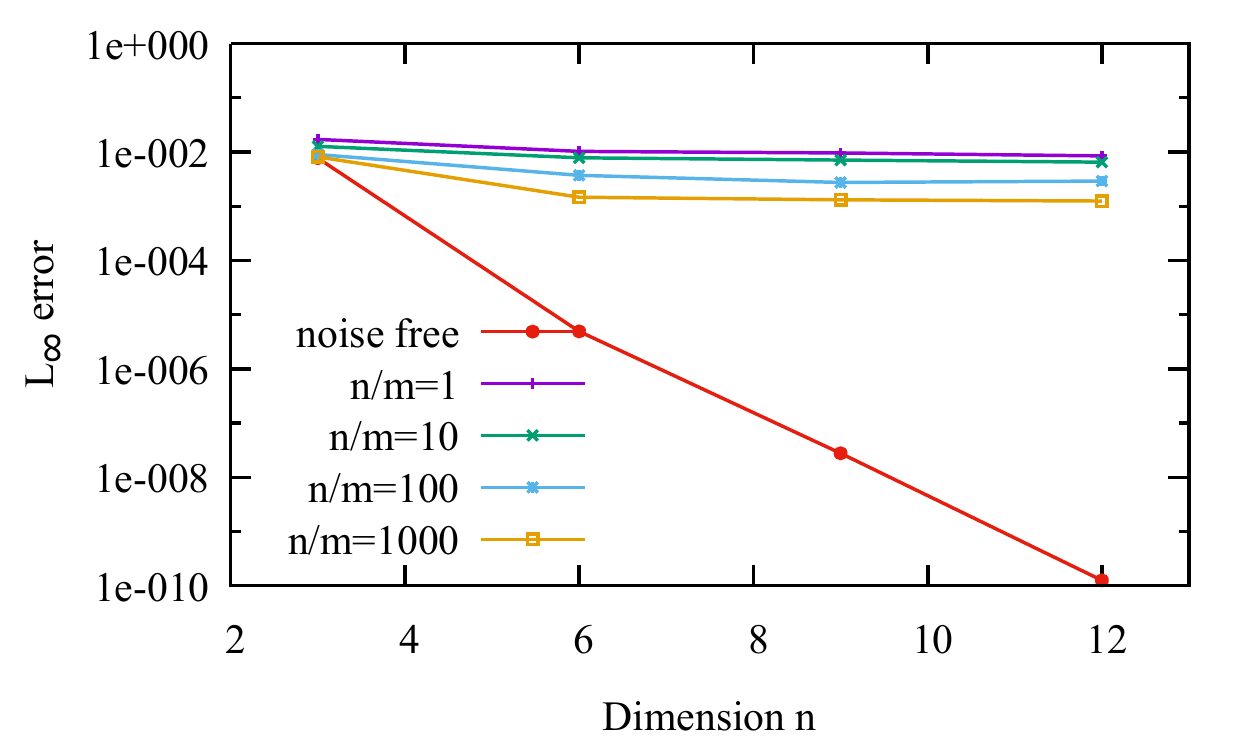}
\caption{CS-GEIM with different $n/m$ ratio, the input noise level is $10^{-2}$, for Case I.}
\label{fig:2d:moremeasurement}
\end{figure}

\section{Conclusions and future works}
\label{sec:generalization}

We have presented some results obtained the Empirical Interpolation Method (EIM) for the reconstruction of the whole field, solution to a simple but representative problem in nuclear reactor physics as an example of a set of parameterized functions. With EIM, a high accuracy
can be got in reconstructing the physical fields, and also a better
sensors deployment is proposed with which most information
can be extracted in a given precision even if only part of the field (either in space or in component) is omitted in the measurement process. Then an improved Empirical Interpolation Method (CS-(G)EIM) is proposed. With CS-(G)EIM, i) the behavior of the interpolant is improved when measurements suffer from noise, ii) the error is dramatically improved in noisy extrapolation case, iii) it is possible to decrease the error by increasing the number of measurements. 

Further works and perspective are ongoing: i)  mathematical analysis of the stable and accurate behavior of this stabilized approach, ii)   in this work, our first assumption is the model is perfect (i.e. we work on in silico solutions, a broader class of methods which couple reduced models with measured data named PBDW \cite{MPPY2015} are able to correct the bias of the model and use real data.

		\bibliographystyle{unsrt}
\bibliography{literature}

\end{document}